\documentclass[12pt]{amsproc}
\usepackage[cp1251]{inputenc}
\usepackage[T2A]{fontenc} 
\usepackage[russian]{babel} 
\usepackage{amsbsy,amsmath,amssymb,amsthm,graphicx,color} 
\usepackage{amscd}

\def\ge{\geqslant}

\newtheorem{prop}{Предложение}

\theoremstyle{definition}

\theoremstyle{remark}

\begin {document}
\centerline{УДК 512.547.2}
\unitlength=1mm
\title[Двойственные характеры алгебры Гекке]
{Двойственные характеры алгебры Гекке и подкрученные следы Маркова}
\author{Г. Г. Ильюта}
\email{ilyuta@mccme.ru}
\address{Московский Государственный Гуманитарный Университет им. М. А. Шолохова}
\thanks{Работа поддержана грантами РФФИ-10-01-00678 и НШ-4850.2012.1.}
\maketitle
\begin{abstract} 
В статье определены подкрученные следы Маркова на алгебрах Гекке конечных групп Кокстера. На характерах алгебры Гекке симметрической группы определена двойственность и доказано, что подкрученные следы Маркова двойственны неприводимым характерам алгебры Гекке. Этот факт приводит к простой формуле для элементов таблицы характеров подкрученных следов Маркова и аналогичному правилу Старки соотношению для таблицы характеров.
\end{abstract}

 Определение следа Маркова обобщено для алгебр Гекке конечных групп Кокстера \cite{6} и некоторых комплексных групп отражений \cite{9}. Для любого неприводимого характера каждой из этих алrебр мы определим подкрученный след Маркова (след из \cite{6} отвечает $q$-аналогу тривиального характера группы). В случае алгебры Гекке $H_n$ симметрической группы $S_n$ мы определим подкрученные следы Маркова также с помощью специализации внутреннего произведения функций Шура и квантового характеристического отображения Фробениуса \cite{4}. На характерах алгебры Гекке $H_n$ мы определим двойственность и докажем, что подкрученные следы Маркова двойственны неприводимым характерам алгебры Гекке $H_n$, в частности, классический след Маркова и $q$-аналог тривиального характера группы $S_n$ двойственны. Неформально говоря, эта двойственность устанавливает своего рода равноправие между симметрической и внешней алгебрами в формулах для весов следа Маркова из \cite{6}, что не очевидно из самих формул. Из утверждения о двойственности вытекает простая формула для элементов таблицы характеров подкрученных следов Маркова, квантовая формула Фробениуса \cite{12} и её супераналог \cite{11} также приводят к формулам для этих элементов, в которых участвуют (супер)симметрические функции Холла-Литтлвуда. Доказано также следующее равенство: произведение матрицы определённых ниже подкрученных рядов Пуанкаре группы $S_n$, таблицы характеров $\chi_q$ алгебры Гекке $H_n$ и некоторой диагональной матрицы равно матрице $\chi_q|_{q=0}$. Это равенство можно считать вариацией правила Старки, связывающего таблицы характеров алгебры Гекке $H_n$ и группы $S_n$ \cite{5} (и его можно вывести из правила Старки и формулы Молина по аналогии с доказательством теоремы 4.1 для серии $A$ в \cite{6}). Возможно, что в такой форме окажется проще обобщить правило Старки для других групп, поскольку формулировка не привязана к симметрической группе (подкрученные ряды Пуанкаре определены для любой конечной группы, а диагональная матрица определяется минимальными длинами элементов в классах сопряжённости группы).

 Классическое характеристическое отображение Фробениуса переводит некоторые специализации подкрученных следов Маркова в $N$-функции Шура $s_\lambda(x^{(N)})$ (в наборе переменных $x^{(N)}$ каждая переменная из набора $x$ повторяется $N$ раз) \cite{1}, \cite{2}, связанные с плетизмом, функциями Джека, парковочными функциями и решёткой нескрещивающихся разбиений. В случае классического следа Маркова такая специализация связана с инвариантами Решетихина-Тураева, приходящими из квантовой группы $U_q(sl(N))$ \cite{7}.

 При изучении специализации следа Маркова -- канонического симметризирующего следа -- обычно использовались не ряды Пуанкаре групп Кокстера, а пропорциональные им многочлены Пуанкаре алгебры коинвариантов (fake degrees) \cite{5}. Мы определим в Примере 6 подкрученные многочлены Пуанкаре алгебры коинвариантов и покажем, что они пропорциональны подкрученным рядам Пуанкаре с тем же коэффициентом пропорциональности.
 
 Мы закончим введение следующим вопросом: можно ли обобщить конструкцию В. Джонса и получить с помощью подкрученных следов Маркова новые инварианты узлов? Одним из возможных путей для поиска ответа на этот вопрос может стать обобщение для подкрученных следов Маркова геометрического определения классического следа Маркова и использование гомологий Хованова-Розанского \cite{13}, \cite{9}.

 1. Определение.

 Для характера $\zeta$ конечной группы обозначим через $M[\zeta]$ матрицу тензорного умножения на $\zeta$, записанную в базисе неприводимых характеров $\zeta_1\equiv 1,\zeta_2,\dots$
$$
\zeta\otimes \zeta_i=M[\zeta]^i_1\zeta_1+M[\zeta]^i_2\zeta_2+\dots
$$
Пусть $\zeta_S^m$ и $\zeta_A^m$ -- $m$-е симметрическая и внешняя степени характера $\zeta$ и 
$$
\zeta_S(q)=\sum_{j\ge 0}\zeta_S^jq^j, \zeta_A(q)=\sum_{j\ge 0}\zeta_A^jq^j.
$$

 В \cite{6} определён след Маркова $\tau$ на алгебре Гекке конечной группы Кокстера. Пусть $\zeta$ -- характер стандартного представления этой группы в $n$-мерном пространстве $V$. В \cite{6} доказана формула 
$$
\tau=\left(\frac{1-q}{1+r}\right)^nPLC_q,
$$
где $C_q$ -- столбец неприводимых характеров алгебры Гекке, $L$ -- матрица введённого Дж. Люстигом преобразования Фурье (для симметрической группы матрица $L$ единичная), $P$ -- строка рядов Пуанкаре тензорного произведения симметрической и внешней алгебр пространства $V$. Из соотношений ортогональности для характеров следует, что строка $P$ совпадает со строкой матрицы $M[\zeta_S(q)\otimes \zeta_A(r)]$, отвечающей единичному характеру. Эти факты мотивируют определение столбца $(\tau_q^{\zeta_i})$ подкрученных следов  Маркова:
$$
(\tau_q^{\zeta_i})=\left(\frac{1-q}{1+r}\right)^nM[\zeta_S(q)\otimes \zeta_A(r)]LC_q.
$$
Результаты \cite{6} обобщены для некоторых комплексных групп отражений \cite{9}, определение подкрученных следов  Маркова годится и для них.

 Порождающие элементы $T_w$ алгебры Гекке $H_n$ индексируются элементами $w$ группы $S_n$, а её неприводимые характеры $\chi^\lambda_q$ (как и неприводимые характеры группы $S_n$) взаимно однозначно соответствуют разбиениям
$$
\lambda=(\lambda_1,\lambda_2, \dots, \lambda_{l(\lambda)}), \lambda_1 \ge \lambda_2 \ge \dots \ge \lambda_{l(\lambda)}>0, \sum \lambda_i=n.
$$
Для разбиения $\lambda$ обозначим через $\lambda'$ сопряжённое разбиение и через $w_\lambda$ элемент Кокстера подгруппы Юнга $S_\lambda=S_{\lambda_1}\times \dots\times S_{\lambda_{l(\lambda)}}$ \cite{5}.

\begin{prop}\label{prop1} Подкрученные следы Маркова $\tau_q^\gamma$, ${\gamma \vdash n}$, на алгебре Гекке $H_n$ симметрической группы $S_n$ связаны с внутренним произведением функций Шура ${s_\lambda*s_\gamma}$ и квантовым характеристическим отображением Фробениуса $ch_q$ формулами
$$
\tau_q^\gamma=\left(\frac{1-q}{1+r}\right)^n \sum_{\lambda \vdash n} s_\lambda*s_\gamma\left(\frac{1+r}{1-q}\right)\chi^\lambda_q
$$
$$
=\left(\frac{1-q}{1+r}\right)^n ch_q^{-1}s_\lambda\left(\frac{1+r}{1-q}y\right).
$$
\end{prop}

 Доказательство. Для симметрической группы строка $P$ совпадает со строкой, состоящей из специализаций функций Шура \cite{6}, \cite{10},I.3.Ex.3,
$$
s_\gamma\left(\frac{1+r}{1-q}\right)=\prod_{(i,j)\in \gamma}\frac{q^i+rq^j} {1-q^{\gamma_i+\gamma'_j-i-j+1}},
$$
произведение берётся по клеткам $(i,j)$ таблицы Юнга разбиения $\lambda$. По определению \cite{10}
$$
s_\lambda*s_\gamma(x)=\sum_{\nu \vdash n}M[\chi^\gamma]_\lambda^\nu s_\nu(x).
$$
Из соотношений ортогональности для характеров вытекают равенства 
$$
M[\zeta]^i_j=(\zeta\otimes \zeta_i,\zeta_j) =((\sum(\zeta,\zeta_k)\zeta_k)\otimes \zeta_i,\zeta_j)    
$$
$$
=\sum(\zeta,\zeta_k)(\zeta_k\otimes \zeta_i,\zeta_j)==\sum M[\zeta]^1_kM[\zeta_i]^k_j.
$$
Поэтому
$$
s_\lambda*s_\gamma\left(\frac{1+r}{1-q}\right)=M[\zeta_S(q)\otimes \zeta_A(r)]^\lambda_\gamma.
$$
Второе равенство получается после применения к обеим частям формулы \cite{12}
$$
s_\gamma(xy)=\sum_{\lambda \vdash n}s_\gamma*s_\lambda(x)s_\lambda(y),
$$
специализации $x\to (1+r)/(1-q)^{-1}$ и отображения $ch_q^{-1}$ по переменным $y$, которое переводит функции Шура в неприводимые характеры алгебры Гекке $H_n$ \cite{4}.

\begin{prop}\label{prop2} Для любого $\gamma \vdash n$ подкрученный след Маркова $\tau_q^\gamma$, характеристическое отображение Фробениуса $ch$ и $N$-функция Шура $s_\gamma(x^{(N)})$ связаны формулой
$$
ch(\lim_{q\to 1}(\tau^\gamma|_{r=-q^N})=s_\gamma(x^{(N)}).
$$
\end{prop}

 Доказательство. При $q=1$ неприводимые характеры алгебры Гекке $H_n$ переходят в неприводимые характеры группы $S_n$ \cite{5}, которые характеристическое отображением Фробениуса переводит в функции Шура \cite{10}. Поэтому утверждение вытекает из формулы для $N$-функций Шура  \cite{1}, \cite{2}
$$
s_\lambda(x^{(N)})=\sum_{\gamma,\nu}\prod_{(i,j)\in \nu}\frac{N+j-i} {\gamma_i+\gamma'_j-i-j+1}M[\chi^\gamma]_\lambda^\nu s_\gamma(x),
$$
и равенства
$$
\lim_{q\to 1}\biggl(\prod_{(i,j)\in \gamma}\frac{q^i+rq^j} {1-q^{\gamma_i+\gamma'_j-i-j+1}}\biggl|_{r=-q^N}\biggr)=\prod_{(i,j)\in \gamma}\frac{N+j-i} {\gamma_i+\gamma'_j-i-j+1}.
$$

 2. Двойственность.

 Для характера $\zeta_q$ алгебры Гекке $H_n$ определим двойственный характер $\hat \zeta_q$ формулой (положим $t=-r$)
$$  
(1-t)^{l(w_\lambda)}\hat \zeta_q(T_{w_\lambda})= (1-q)^{l(w_\lambda)}\zeta_t(T_{w_\lambda}),\lambda\vdash n. 
$$
Любой характер определяется своими значениями на элементах $T_{w_\lambda}$ \cite{12}. По определению $\zeta_t(T_{w_\lambda})=\zeta_q(T_{w_\lambda})|_{q=t}$ и аналогично определяется $ch_t$.

\begin{prop}\label{prop3} Образы двойственных характеров при квантовом характеристическом отображении Фробениуса $ch_q$ связаны формулой
$$
(1-t)^n ch_q(\hat \zeta_q)((q-1)x)=(1-q)^n ch_t(\zeta_t)((t-1)x).
$$
\end{prop}

 Доказательство. Используя разложение образа по мономиальным симметрическим функциям $m_\lambda(x)$ \cite{4}
$$
ch_q(\zeta_q)((q-1)x)=\sum_{\lambda \vdash n}(q-1)^{l(\lambda)} \zeta_q(T_{w_\lambda})m_\lambda(x).
$$
и равенсто $l(\lambda)+l(w_\lambda)=n$, получим
$$
(1-t)^n ch_q(\hat \zeta_q)((q-1)x)=(1-t)^n\sum_{\lambda \vdash n}(q-1)^{l(\lambda)} \hat \zeta_q(T_{w_\lambda})m_\lambda(x)
$$
$$
=(1-q)^n\sum_{\lambda \vdash n}(t-1)^{l(\lambda)} \zeta_t(T_{w_\lambda})m_\lambda(x)=(1-q)^n ch_t(\zeta_t)((t-1)x).
$$

\begin{prop}\label{prop4} Для любого $\gamma \vdash n$ характеры $\chi_q^\gamma$ и $\tau_q^\gamma$ двойственны, тем самым
$$
(1-t)^{l(w_\beta)}\tau_q^\lambda(T_{w_\beta})= (1-q)^{l(w_\beta)}\chi_t^\lambda(T_{w_\beta}).
$$
\end{prop}

 Доказательство. Используя формулу из \cite{4} и Предложение 1, имеем
$$
\sum_{\beta \vdash n} (1-t)^n(q-1)^{l(\beta)} \tau_q^\lambda(T_{w_\beta})m_\beta(x)=(1-t)^n ch_q(\tau_q^\lambda)((q-1)x)
$$
$$
=(1-t)^n\left(\frac{1-q}{1+r}\right)^n s_\lambda\left((q-1)\frac{1+r}{1-q}x\right)
$$
$$
=(1-q)^n s_\lambda((t-1)x)=\sum_{\beta \vdash n} (1-q)^n (t-1)^{l(\lambda)}\chi_q^\lambda(T_{w_\beta})m_\beta(x).
$$

\begin{prop}\label{prop5} Условие двойственности подкрученных следов Маркова и неприводимых характеров алгебры Гекке $H_n$ равносильно следующему соотношению для таблицы характеров $\chi_q=(\chi_q^\gamma (T_{w_\beta}))$ алгебры Гекке $H_n$
$$
M[\zeta_S(q)]\chi_q D[(1-q)^{l(\lambda)}]=\chi_q|_{q=0},
$$
где $D[(1-q)^{l(\lambda)}]$ -- диагональная матрица с элементами $(1-q)^{l(\lambda)}$, $\lambda \vdash n$, на диагонали.
\end{prop}

 Доказательство. Поскольку $\zeta_S(q)\otimes \zeta_A(-q)\equiv 1$ \cite{3}, то
$$
M[\zeta_S(q)\otimes \zeta_A(r)]=M[\zeta_S(q)]M[\zeta_S(-r)]^{-1}.
$$
Известно, что все матрицы тензорного умножения на характеры диагонализируются в одном и том же базисе и потому коммутируют. Записывая условие двойственности в матричной форме, имеем
$$
\chi_t D[(1-q)^{l(w_\lambda)}]D[(1-t)^{l(w_\lambda)}]^{-1}= (\tau_q^\gamma(T_{w_\beta}))
$$
$$
=\left(\frac{1-q}{1+r}\right)^n M[\zeta_S(q)]M[\zeta_S(-r)]^{-1}\chi_q.
$$
$$
(1-q)^n M[\zeta_S(q)]\chi_q D[(1-q)^{l(w_\lambda)}]^{-1}=
(1-t)^n M[\zeta_S(t)]\chi_t D[(1-t)^{l(w_\lambda)}]^{-1}.
$$
Поэтому матрица $M[\zeta_S(q)]\chi_q D[(1-q)^{l(\lambda)}]$ не зависит от $q$, а значит
$$
M[\zeta_S(q)]\chi_q D[(1-q)^{l(\lambda)}]=(M[\zeta_S(q)]\chi_q D[(1-q)^{l(\lambda)}])|_{q=0}=\chi_q|_{q=0}.
$$
Ясно, что это рассуждение обратимо.

 3. Примеры.

 Пример 1. Для одномерных характеров $\chi_q^{(n)}$, $\chi_q^{(1^n)}$ алгебры Гекке $H_n$ имеем
$$
\chi_q^{(n)}(T_w)=q^{l(w)},\chi_q^{(1^n)}(T_w)=(-1)^{l(w)}
$$
Поэтому
$$
\tau^{(n)}(T_{w_\beta})=z^{l(w_\beta)},\tau^{(1^n)}(T_{w_\beta})=\left(\frac{z}{r}\right)^{l(w_\beta)},
$$
где $z(1+r)=(q-1)r$.

 Пример 2. Согласно \cite{5}, Prop. 9.4.1.,
$$
\chi_q^{\lambda'}(T_w)=(-q)^{l(w)}\chi_q^\lambda(T_w)|_{q\to q^{-1}}.
$$
Поэтому
$$
\tau_z^{\lambda'}(T_{w_\beta})=(-z)^{l(w_\beta)}\chi_t^\lambda(T_{w_\beta})|_{t\to t^{-1}}.
$$

 Пример 3. Согласно \cite{5}, Lemma 9.2.1., для базисных транспозиций $s_i=(ii+1)\in S_n$ и любого разбиения $\lambda$ имеем равенство
$$
\chi_q^\lambda(T_{s_i})=\frac{q}{2}(\chi^\lambda(1)+\chi^\lambda(s_i))-\frac{1}{2}(\chi^\lambda(1)-\chi^\lambda(s_i)).
$$
Поскольку $w_{(21^{n-2})}=s_1$, то
$$
\tau_z^\lambda(T_{w_{(21^{n-2})}})=\frac{z}{2}(\chi^\lambda(1)+\chi^\lambda(s_i))-\frac{z}{2t}(\chi^\lambda(1)-\chi^\lambda(s_i)).
$$

 Пример 4. Симметрические функции Холла-Литтлвуда $q_\lambda(x;q)=\prod q_{\lambda_i}(x;q)$ определяются производящей функцией 
$$
\sum_{j=0}^\infty q_j(x;q)u^j=\prod_i\frac{1-qx_iu}{1-x_iu}.
$$
Рассмотрим внутреннее произведение обеих частей квантовой формулы Фробениуса \cite{12}
$$
\frac{q^n}{(q-1)^{l(\mu)}}q_\mu(x;q^{-1})=\sum_{\lambda \vdash n}\chi_q^\lambda (T_{w_\mu})s_\lambda(x)
$$
и функции Шура $s_\gamma(x)$ и затем специализируем переменные $x\to (1+r)/(1-q)^{-1}$. Справа получим $\tau_q^\gamma (T_{w_\mu})$. Мы представили элементы таблицы характеров подкрученных следов Маркова как специализации внутренних произведений функций Холла-Литтлвуда и функций Шура.

 Пример 5. Рассмотрим симметрические функции в квантовой формуле Фробениуса как функции от двух наборов переменных $x$ и $y$. Применим к обеим частям по переменным $y$ инволюцию, переставляющую элементарные и полные симметрические функции \cite{10}. Справа получим суперсимметрические функции Шура $s_\lambda(x/y)$ \cite{10},I.5.Ex.23, а слева -- суперсимметрические функции Холла-Литтлвуда $q_\lambda(x/y;q)$ \cite{11}, поскольку из определения инволюции следует, что образом зависящей от переменных $x$ и $y$ производящей функции из Примера 4 будет производящая функция для суперсимметрических функций Холла-Литтлвуда
$$
\sum_{j=0}^\infty q_j(x/y;q)u^j=\prod_i\frac{1-qx_iu}{1-x_iu}\prod_j\frac{1-y_ju}{1-qy_ju}.
$$
Мы получили супераналог квантовой формулы Фробениуса (в \cite{11} приводится доказательство, обобщающее доказательство из \cite{12} и использующее супераналог двойственности Шура-Вейля). Отсюда аналогично Примеру 4 получаем формулу для элементов таблицы характеров подкрученных следов Маркова. Специализация $x\to (1+r)/(1-q)^{-1}$ в контексте суперсимметрических функций изучалась в \cite{8}.

 Пример 6. Пусть $d_i$ -- размерность пространства инвариантных многочленов степени $i$ на пространстве $V$. Алгебра коинвариантов получается факторизацией симметрической алгебры пространства $V$ по инвариантным многочленам положительной степени, по аналогии с $\zeta_S(q)$ и $\zeta_A(q)$ обозначим через 
$$
\zeta_C(q)=\sum_{j\ge 0}\zeta_C^jq^j
$$
градуированный характер на алгебре коинвариантов. Поскольку
$$
\sum_{i\ge 0}d_iq^i=\prod_j\frac{1}{1-q^{e_j}},
$$
где $e_1,e_2,\dots$ -- показатели группы Кокстера, и, согласно \cite{3}, Prop. 11.1.1.,
$$
\zeta_S^m=\sum_{i+j=m}d_i\zeta_C^j,
$$
то
$$
(\sum_{i\ge 0}\zeta_S^iq^i\otimes \zeta_k,\zeta_l)=\sum_{i\ge 0}d_iq^i(\sum_{j\ge 0}\zeta_C^jq^j \otimes \zeta_k,\zeta_l),
$$
$$
M[\zeta_C(q)]=\prod (1-q^{e_j})M[\zeta_S(q)].
$$

\bigskip

\end {document}